\documentclass[leqno,12pt,twoside]{amsart}
\usepackage{geometry}

\usepackage{times}

\usepackage[all]{xy} 
\usepackage{amsmath, amssymb, amsfonts, latexsym, mdwlist, amsthm}
\usepackage{subfig}
\usepackage{graphicx}
\usepackage{wrapfig}

\usepackage{hyperref}

\usepackage{tikz}
\usetikzlibrary{calc,trees,positioning,arrows,chains,shapes.geometric,%
    decorations.pathreplacing,decorations.pathmorphing,shapes,%
    matrix,shapes.symbols}

\tikzset{
>=stealth',
  punktchain/.style={
    rectangle,
    rounded corners,
    draw=black, thick,
    minimum height=3em,
    text centered,
    on chain},
  line/.style={draw, thick, <-},
  element/.style={
    tape,
    top color=white,
    bottom color=blue!50!black!60!,
    minimum width=8em,
    draw=blue!40!black!90, very thick,
    text width=10em,
    minimum height=3.5em,
    text centered,
    on chain},
  every join/.style={->, thick,shorten >=1pt},
  decoration={brace},
  tuborg/.style={decorate},
  tubnode/.style={midway, right=2pt},
}

\usepackage{paralist}
\setdefaultenum{(a)}{(i)}{}{}
\usepackage{enumitem} 

\def\C{\ensuremath{\mathbb{C}}}
\def\D{\ensuremath{\mathbb{D}}}
\def\H{\ensuremath{\mathbb{H}}}

\def\P{\ensuremath{\mathbb{P}}}
\def\Q{\ensuremath{\mathbb{Q}}}
\def\R{\ensuremath{\mathbb{R}}}

\def\ch{\mathop{\mathrm{ch}}\nolimits}

\def\Coh{\mathop{\mathrm{Coh}}\nolimits}

\def\Hom{\mathop{\mathrm{Hom}}\nolimits}

\def\RlHom{\mathop{\mathbf{R}\mathcal Hom}\nolimits}

\def\Num{\mathop{\mathrm{Num}}\nolimits}
\def\NS{\mathop{\mathrm{NS}}\nolimits}

\newenvironment{Prf}{\textit{Proof.}\/}{\hfill$\Box$}

\def\MG13{\ensuremath{{\mathcal M}_{\Gamma_1(3)}}}
\def\tildeMG13{\ensuremath{\widetilde{\mathcal M}_{\Gamma_1(3)}}}

\def\into{\ensuremath{\hookrightarrow}}
\def\onto{\ensuremath{\twoheadrightarrow}}

\def\blank{\underline{\hphantom{A}}}

\def\Db{\mathrm{D}^{\mathrm{b}}}

\newcommand\stv[2]{\left\{#1\,\colon\,#2\right\}}

\makeatletter
\newtheorem*{rep@theorem}{\rep@title}
\newcommand{\newreptheorem}[2]{%
\newenvironment{rep#1}[1]{%
 \def\rep@title{#2 \ref{##1}}%
 \begin{rep@theorem}}%
 {\end{rep@theorem}}}
\makeatother

\newtheorem{Thm}{Theorem}[section]
\newreptheorem{Thm}{Theorem}
\newtheorem{Prop}[Thm]{Proposition}

\newtheorem{Lem}[Thm]{Lemma}
\newtheorem{Cor}[Thm]{Corollary}
\newreptheorem{Cor}{Corollary}
\newtheorem{Con}[Thm]{Conjecture}
\newreptheorem{Con}{Conjecture}

\newtheorem{Claim}{Claim}

\newtheorem{thm-int}{Theorem}

\theoremstyle{definition}
\newtheorem{Def-s}[Thm]{Definition}
\newtheorem{Def}[Thm]{Definition}

\newtheorem{Ex}[Thm]{Example}

\def\C{\ensuremath{\mathbb{C}}}
\def\D{\ensuremath{\mathbb{D}}}
\def\H{\ensuremath{\mathbb{H}}}

\def\P{\ensuremath{\mathbb{P}}}
\def\Q{\ensuremath{\mathbb{Q}}}
\def\R{\ensuremath{\mathbb{R}}}

\def\AA{\ensuremath{\mathcal A}}
\def\BB{\ensuremath{\mathcal B}}
\def\CC{\ensuremath{\mathcal C}}

\def\EE{\ensuremath{\mathcal E}}
\def\FF{\ensuremath{\mathcal F}}
\def\GG{\ensuremath{\mathcal G}}
\def\HH{\ensuremath{\mathcal H}}
\def\II{\ensuremath{\mathcal I}}

\def\OO{\ensuremath{\mathcal O}}

\def\QQ{\ensuremath{\mathcal Q}}
\def\TT{\ensuremath{\mathcal T}}

\def\cht{\ensuremath{\ch^B}}

\def\EEE{\mathfrak E}

\newcommand{\ignore}[1]{}

\begin{document}



\title[Bogomolov-Gieseker inequality for the projective space]{A generalized Bogomolov-Gieseker inequality for the three-dimensional projective space}

\author{Emanuele Macr\`i}
\address{Department of Mathematics, The Ohio State University, 231 W 18th Avenue, Columbus, OH 43210-1174, USA}
\email{macri.6@math.osu.edu}
\urladdr{http://www.math.osu.edu/~macri.6/}



\begin{abstract}
A generalized Bogomolov-Gieseker inequality for tilt-stable complexes on a smooth projective threefold was conjectured by Bayer, Toda, and the author.
We show that such inequality holds true in general, if it holds true when the polarization is sufficiently small.
As an application, we prove it for the three-dimensional projective space.
\end{abstract}

\maketitle

\setcounter{tocdepth}{1}

\section{Introduction}\label{sec:intro}

The notion of tilt-stability, for objects in the derived category of a smooth projective threefold, was introduced in \cite{BMT:3folds-BG}, based on \cite{Bridgeland:K3, Aaron-Daniele}.
In \cite[Conjecture 1.3.1]{BMT:3folds-BG} (Conjecture \ref{con:strongBG} of the present paper), we proposed a generalized Bogomolov-Gieseker inequality (BG inequality, for short) for tilt-stable objects.
The main application for tilt-stability was to have an auxiliary notion of stability to construct Bridgeland stability conditions.
The generalized BG inequality is precisely the missing ingredient to being able to show the existence of Bridgeland stability conditions.

In this note, we prove such inequality in the case of the projective space $\P^3$.

\begin{Thm}\label{thm:main}
The generalized Bogomolov-Gieseker inequality for tilt-stable objects in $\Db(\P^3)$ holds.
\end{Thm}

This gives the first example when the generalized BG inequality is proved in full generality.
As a corollary, by \cite{BMT:3folds-BG}, we can also describe a large open subset of the space of stability conditions on $\Db(\P^3)$.
It would be very interesting to study how moduli spaces of Bridgeland semistable objects vary when varying the stability condition (very much like the situation described in \cite{ABCH:MMP, MaciociaMeachan, LoQin:miniwalls, MYY2, YY:abeliansurfaces, BM:Projectivity,Toda:MMPxSurfaces,Yoshioka:GeneralizedKummer}, for the case of surfaces).
The behavior at the ``large volume limit point'' is described in \cite[Section 6]{BMT:3folds-BG}.

The idea of the proof of Theorem \ref{thm:main} goes as follows.
For a smooth projective threefold $X$, the notion of tilt-stability depends on two parameters, namely two divisor classes $B,\omega\in\mathrm{NS}_{\R}(X)$, with $\omega$ ample.
In this paper we prove a general result, Proposition \ref{prop:reduction}: To show the generalized BG inequality for all $B$ and $\omega$ can always be reduced to showing it for $\omega$ ``arbitrarily small'', uniformly in $B$.

For $X=\P^3$, the case in which $\omega$ is small was essentially proved in \cite[Theorem 8.2.1]{BMT:3folds-BG}.
More precisely, for simplicity, in \cite{BMT:3folds-BG} only the case $B=0$ was considered.
Proposition \ref{prop:P3Bound} generalizes that argument to arbitrary $B$.
Together with Proposition \ref{prop:reduction}, this completes the proof of Theorem \ref{thm:main}.

The interest for a general proof of the generalized BG inequality, besides for the existence of Bridgeland stability conditions, relies on its consequences.
Indeed, if we assume such inequality to be true, we would have:
\begin{itemize}
\item A proof of Fujita's Conjecture for threefolds, \cite{BBMT:Fujita};
\item A mathematical formulation of Denef-Moore's formula derived in the study of Ooguri-Strominger-Vafa's conjecture, relating black hole entropy and topological string, \cite{Toda:DenefMoore};
\item The possibility to realize extremal contractions for threefolds as moduli spaces of semistable objects in the derived category, \cite{Toda:ExtremalContractions}.
\end{itemize}
We also mention that in the paper \cite{Polishchuk:LIobjects} the existence of Bridgeland stability conditions on abelian threefolds is tested on a class of objects (called Lagrangian-Invariant objects).

Finally, in \cite{BMT:3folds-BG} it was pointed out a strict relation between the generalized BG inequality and Castelnuovo's inequality for curves in $\P^3$.
In Section \ref{sec:applications} of this paper, we show that Theorem \ref{thm:main} gives, as an immediate corollary, a weaker version of Castelnuovo's theorem \cite[IV, 6.4]{Hartshorne}.

A survey on Bridgeland stability conditions and further problems and applications can be found in \cite{Bridgeland:spaces, stability-tour,Huybrechts:Survey, Toda:Survey}.

\subsection*{Notation}

In this paper, we will always denote by $X$ a smooth projective threefold over the complex numbers and by $\Db(X)$ its bounded derived category of coherent sheaves.
The Chow groups of $X$ modulo numerical equivalence are denoted by $\mathrm{Num}(X)$.
In particular, the N\'eron-Severi group $\mathrm{NS}(X)=\mathrm{Num}^1(X)$.

For an abelian group $G$ and a field $k(=\Q,\R,\C)$, we denote by $G_k$ the $k$-vector space $G\otimes k$.

\subsection*{Acknowledgements}

I would like to thank Arend Bayer and Ciaran Meachan for very useful discussions and Roman Bezrukavnikov for pointing out the idea of using Theorem \ref{thm:main} to show a version of Castelnuovo's Theorem.
This paper was completed during my stay at the Mathematical Institute of the University of Bonn whose warm hospitality is gratefully acknowledged.
I was partially supported by the NSF grant DMS-1160466 and, during the visit to Bonn, by the grant SFB/TR 45.

\section{The reduction argument}\label{sec:reduction}

In this section we give a brief recall on the notion of tilt stability, following \cite{BMT:3folds-BG}.
We show how to reduce the proof of \cite[Conjecture 1.3.1]{BMT:3folds-BG}, when $\omega$ and $B$ are ``parallel'', to the case in which the polarization is ``sufficiently small''.

\subsection{Tilt stability}\label{subsec:TiltStability}

Let $X$ be a smooth projective threefold over $\C$, and let $H\in\mathrm{NS}(X)$ be an ample divisor class.
For a pair
\begin{align*}
&\omega = \alpha \cdot H, \qquad \alpha\in\R_{>0}\\
&B = \beta \cdot H, \qquad \beta\in\R,
\end{align*}
we define a slope function $\mu_{\omega, B}$ for coherent sheaves on $X$ in the usual way:
For $E\in\Coh(X)$, we set
\[
\mu_{\omega, B}(E) = \begin{cases} 
 + \infty, & \text{ if }\cht_0(E)=0,\\
\ & \\
\frac{\omega^2 \cht_1(E)}{\omega^3 \cht_0(E)}, & \text{ otherwise,}
 \end{cases}
\]
where $\cht(E) = e^{-B} \ch(E)$ denotes the Chern character twisted by $B$.
Explicitly:
\begin{align*}
&\cht_0=\ch_0 & \cht_2=\ch_2-B\ch_1+\frac{B^2}{2}\ch_0 & \\
&\cht_1=\ch_1-B\ch_0 &\cht_3=\ch_3-B\ch_2+\frac{B^2}{2}\ch_1-\frac{B^3}{6}\ch_0 &.
\end{align*}

A coherent sheaf $E$ is slope-(semi)stable (or $\mu_{\omega,B}$-(semi)stable) if, for all subsheaves $F \into E$, we have
\[
\mu_{\omega, B}(F) < (\le) \mu_{\omega, B}(E/F).
\]

Due to the existence of Harder-Narasimhan filtrations (HN-filtrations, for short) with respect to slope-stability,
there exists a \emph{torsion pair} $(\TT_{\omega, B}, \FF_{\omega, B})$ defined 
as follows:
\begin{align*}
\TT_{\omega, B} &= \stv{E \in \Coh X}
{\text{any quotient $E \onto G$ satisfies $\mu_{\omega, B}(G) > 0$}} \\
\FF_{\omega, B} &= \stv{E \in \Coh X}
{\text{any subsheaf $F \into E$ satisfies $\mu_{\omega, B}(F) \le 0$}}
\end{align*}
Equivalently, $\TT_{\omega, B}$ and $\FF_{\omega, B}$ are the extension-closed
subcategories of $\Coh X$ generated by slope-stable sheaves of positive or non-positive slope,
respectively.

\begin{Def}\label{def:BB}
We let $\Coh^{\omega, B}(X) \subset \Db(X)$ be the extension-closure
\[ \Coh^{\omega, B}(X) = \langle \TT_{\omega, B}, \FF_{\omega, B}[1] \rangle.
\]
\end{Def}

The category $\Coh^{\omega,B}(X)$ depends only on $\omega$ via $H$.
Hence, to simplify notation, since for us $B$ is also a multiple of $H$, we denote it by $\Coh^B(X)$.
By the general theory of torsion pairs and tilting \cite{Happel-al:tilting},
$\Coh^{B}(X)$ is the heart of a bounded t-structure on $\Db(X)$.

By using the classical Bogomolov-Gieseker inequality and Hodge Index theorem, we
can define the following slope function on $\Coh^{B}(X)$:
For $E\in\Coh^{B}(X)$, we set
\[
\nu_{\omega, B}(E) = \begin{cases} 
 + \infty, & \text{ if }\omega^2 \cht_1(E) = 0,\\
\ & \\
\frac{\omega \cht_2(E) - \frac 12 \omega^3 \cht_0(E)}{\omega^2 \cht_1(E)}, & \text{ otherwise.}
\end{cases}
\]

\begin{Def} \label{def:tilt-stable}
An object $E \in \Coh^{B}(X)$ is \emph{tilt-(semi)stable} if, for all non-trivial subobjects
$F \into E$, we have
\[
\nu_{\omega, B}(F) < (\le) \nu_{\omega, B}(E/F).
\]
\end{Def}

The following is our main conjecture:

\begin{Con}[{\cite[Conjecture 1.3.1]{BMT:3folds-BG}}] \label{con:strongBG}
For any $\nu_{\omega, B}$-semistable object $E\in \Coh^{B}(X)$ satisfying 
$\nu_{\omega, B}(E) = 0$, we have the following inequality
\begin{equation} \label{eq:strongBG}
\cht_3(E) \le \frac{\omega^2}{6}\cht_1(E).
\end{equation}
\end{Con}

The original definition of tilt-stability in \cite{BMT:3folds-BG} was given when $\alpha,\beta\in\Q$ (actually it was slightly more general, allowing $\omega$ and $B$ to be arbitrary, and $\omega$ had a different parameterization $\omega \mapsto \sqrt{3}\cdot \omega$).
The extension to $\R$ is the content of the following proposition, which we recall for later use:

\begin{Prop}[{\cite[Corollary 3.3.3]{BMT:3folds-BG}}] \label{prop:openness}
Let $\mathrm{St} \subset \NS_\R(X) \times \NS_\R(X)$ be the subset of pairs of real 
classes $(\omega, B)$ for which $\omega$ is ample.
There exists a notion of ``tilt-stability'' for every $(\omega, B) \in \mathrm{St}$.
For every object $E$, the set of $(\omega, B)$ for which $E$ is
$\nu_{\omega, B}$-stable defines an open subset of $\mathrm{St}$.
\end{Prop}

\begin{Def}\label{def:discriminant}
We define the \emph{generalized discriminant}
\[
\overline{\Delta}_H := (H^2 \cht_1)^2 -2H^3 \cht_0 \cdot (H\cht_2).
\]
\end{Def}

The generalized discriminant is independent of $\beta$.
Indeed, by expanding the definition, we have
\begin{equation*}
\begin{split}
\overline{\Delta}_H & = (H^2(\ch_1 - \beta \ch_0 H))^2 - 2 H^3 \ch_0 \cdot H(\ch_2 - \beta H\ch_1 + \frac{\beta^2}{2}\ch_0H^2)\\
     		& = (H^2\ch_1)^2 - 2 (H^2\ch_1) H^3 \beta \ch_0 + \beta^2 (\ch_0)^2 (H^3)^2 - 2 H^3 \ch_0 (H\ch_2)\\
		& \quad + 2 (H^2\ch_1) H^3 \beta \ch_0 - \beta^2 (\ch_0)^2(H^3)^2\\
		& = (H^2\ch_1)^2 - 2 H^3 \ch_0 (H\ch_2).
\end{split}
\end{equation*}

The following result will be the key ingredient in our proof:

\begin{Thm}[{\cite[Corollary 7.3.2]{BMT:3folds-BG}}]\label{thm:BogNoCh2}
For any $\nu_{\omega, B}$-semistable object $E\in \Coh^{B}(X)$, we have
\[
\overline{\Delta}_{H}(E)\geq 0.
\]
\end{Thm}

\subsection{Reduction to small $\omega$}\label{subsec:bound}

In this section we prove our reduction result.
We keep the same notation as before, e.g, $\omega=\alpha H$ and $B = \beta H$.
To simplify, we will denote $\nu_{\alpha,\beta}$ for $\nu_{\omega,B}$, $\Coh^{\beta}(X)$, and so on.

\begin{Prop}\label{prop:reduction}
Assume there exists $\overline{\alpha}\in\R_{>0}$ such that, for all $\alpha<\overline{\alpha}$, and for all $\beta\in\R$, Conjecture \ref{con:strongBG} holds.
Then Conjecture \ref{con:strongBG} holds for all $\alpha\in\R_{>0}$ and for all $\beta\in\R$.
\end{Prop}

To prove Proposition \ref{prop:reduction}, we need first to introduce a bit more of notation.
We denote by $\mathbb{H}$ the upper half-plane
\[
\mathbb{H} := \left\{(\beta,\alpha)\in\R^2\,:\,\alpha>0 \right\}.
\]
For a vector $v:=(\ch_0,\ch_1,\ch_2,\ch_3)\in \mathrm{Num}_\Q (X)$, such that $H^2\ch^{\beta}_1>0$, the equation $\nu_{\alpha,\beta}(v)=0$ defines a curve $\CC_v$ in $\mathbb{H}$.
Explicitly, we have
\[
\CC_v:\quad H\ch_2 - \beta (H^2\ch_1) + \frac{\beta^2}{2} H^3 \ch_0 - \frac{\alpha^2}{2} H^3 \ch_0 = 0,
\]
together with the inequality
\[
\beta H^3\ch_0 < H^2\ch_1.
\]
We can divide in two cases:
\begin{align}
& \ch_0=0 \quad \rightsquigarrow \quad \beta = \frac{H\ch_2}{H^2\ch_1},\label{eq:rank0}\\
& \ch_0\neq0 \quad \rightsquigarrow \quad \left(\beta - \frac{H^2\ch_1}{H^3\ch_0} \right)^2 - \alpha^2 = \frac{\overline{\Delta}_H}{(H^3\ch_0)^2}.\label{eq:general}
\end{align}
Hence, if $\overline{\Delta}_H\geq0$, then the tangent line at a point $(\beta_0,\alpha_0)\in\CC$ intersects the line $\alpha=0$ with an angle $\frac{\pi}{4}\leq\theta\leq\frac{\pi}{2}$.

Finally, on the curve $\CC_v$ we can write the inequality \eqref{eq:strongBG} as follows:
\begin{equation}\label{eq:BGonC}
\begin{split}
& \ch_0=0 \quad \rightsquigarrow \quad \ch_3 - \frac{(H\ch_2)^2}{2(H^2\ch_1)}\leq \alpha^2\frac{H^2\ch_1}{6},\\
& \ch_0\neq0 \quad \rightsquigarrow \quad \beta \frac{\overline{\Delta}_H}{H^3 \ch_0} \leq \frac{(H\ch_2)(H^2\ch_1)}{H^3\ch_0} - 3\ch_3.
\end{split}
\end{equation}

\begin{Prf} (Proposition \ref{prop:reduction})
We argue by contradiction.
Assume that there exist $\alpha_0\geq\overline{\alpha}$, $\beta_0\in\R$, an object $E_0\in\Coh^{\beta_0}(X)$ which is $\nu_{\alpha_0,\beta_0}$-stable, such that $\nu_{\alpha_0,\beta_0}(E_0)=0$ and which does not satisfy the inequality in Conjecture \ref{con:strongBG}.

\begin{Claim}\label{claim1}
There exist a sequence $(\beta_n,\alpha_n)\in\mathbb{H}$ and a sequence objects $\{E_n\}_{n\geq0}$ such that
\begin{itemize}
\item $E_n\in\Coh^{\beta_n}(X)\cap\Coh^{\beta_{n+1}}(X)$ is $\nu_{\alpha_n,\beta_n}$-stable,
\item $\nu_{\alpha_n,\beta_n}(E_n)=\nu_{\alpha_{n+1},\beta_{n+1}}(E_n)=0$,
\item $0<H^2\ch^{\beta_{n+1}H}(E_{n+1})< H^2\ch^{\beta_{n+1}H}_1(E_n)$.
\item $E_n$ does not satisfy the inequality \eqref{eq:strongBG},
\item $\alpha_0>\alpha_1>\ldots>\alpha_n>\ldots>0$,
\item $|\beta_{n+1}| \leq |\beta_0| + \alpha_0$.
\end{itemize}
\end{Claim}

\begin{Prf} (Claim \ref{claim1})
We proceed by induction, the case $n=0$ being our assumption.
Assume that we have constructed $E_n$ with the wanted properties.
By Proposition \ref{prop:openness}, the locus in $\mathbb{H}$ where $E_n$ is $\nu_{\alpha,\beta}$-stable is open.
Consider the curve $\CC:=\CC_{\ch(E_n)}\subset\mathbb{H}$ and consider the set $U:=\left\{(\beta,\alpha)\in\CC\,:\, \alpha<\alpha_n\right\}$.
By \eqref{eq:BGonC}, for all $(\beta,\alpha)\in U$, the inequality \eqref{eq:strongBG} is not satisfied for $E_n$.
Since Conjecture \ref{con:strongBG} holds when $\alpha<\overline{\alpha}$, there must exist $(\beta_{n+1},\alpha_{n+1})\in U$ such that $E_n$ is $\nu_{\alpha_{n+1},\beta_{n+1}}$-semistable and it is not $\nu_{\alpha,\beta}$-semistable, for all $(\beta,\alpha)\in U$, with $\alpha<\alpha_{n+1}$.
By looking at the $\nu_{\alpha_{n+1},\beta_{n+1}}$-stable factors of $E_n$ (by \cite[Proposition 5.2.2]{BMT:3folds-BG}, this makes sense in the category $\Coh^{\beta_{n+1}}(X)$), given the additivity of the Chern character, there exists an object $E_{n+1}\in\Coh^{\beta_{n+1}}(X)$ which is $\nu_{\alpha_{n+1},\beta_{n+1}}$-stable, such that $\nu_{\alpha_{n+1},\beta_{n+1}}(E_{n+1})=0$ and which does not satisfy the inequality \eqref{eq:strongBG}.

The final inequality, $|\beta_{n+1}| \leq |\beta_0| + \alpha_0$, follows simply by the fact, observed before, that the tangent line at any point in $\CC$ intersects the line $\alpha=0$ with an angle $\frac{\pi}{4}\leq\theta\leq\frac{\pi}{2}$.
See Figure \ref{fig:SequenceEn}.
\end{Prf}

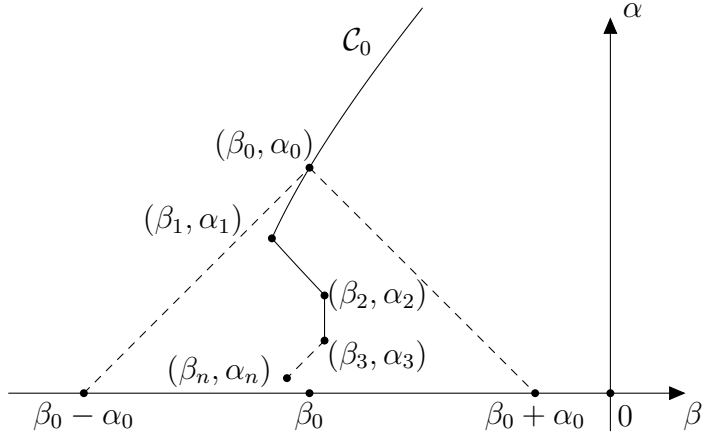
\begin{figure}
\begin{center}

\definecolor{uququq}{rgb}{0.25,0.25,0.25}
\begin{tikzpicture}[line cap=round,line join=round,>=triangle 45,x=1.0cm,y=1.0cm]
\draw[->,color=black] (-8,0) -- (1,0);
\draw[->,color=black] (0,-0.5) -- (0,5);
\fill [color=black] (0,0) circle (1.5pt);
\draw[color=black] (0.2,-0.3) node {$0$};
\fill [color=black] (-1,0) circle (1.5pt);
\draw[color=black] (-1,-0.3) node {$\beta_0+\alpha_0$};
\fill [color=black] (-4,0) circle (1.5pt);
\draw[color=black] (-4,-0.3) node {$\beta_0$};
\fill [color=black] (-7,0) circle (1.5pt);
\draw[color=black] (-7,-0.3) node {$\beta_0-\alpha_0$};
\fill [color=black] (-4,3) circle (1.5pt);
\draw[color=black] (-4.6,3.3) node {$(\beta_0,\alpha_0)$};
\draw[color=black] (1.1,-0.3) node {$\beta$};
\draw[color=black] (0.3,5.1) node {$\alpha$};
\draw[dashed,color=black] (-7,0) -- (-4,3);
\draw[dashed,color=black] (-4,3) -- (-1,0);
\draw (-3.7,5) node[anchor=north west] {$\CC_0$};
\draw [samples=50,domain=-5.5:-3.5,xshift=0cm,yshift=0cm] plot
({-8-\x},{(-16+(\x-1)^2)^(0.5)});
\fill [color=black] (-4.5,2.06) circle (1.5pt);
\draw[color=black] (-5.55,2.3) node {$(\beta_1,\alpha_1)$};
\draw [samples=50,domain=-4.5:-3.8,xshift=0cm,yshift=0cm] plot
({\x},{4.12-(-1+(\x+6.79)^2)^(0.5)});
\fill [color=black] (-3.8,1.3) circle (1.5pt);
\draw[color=black] (-3.1,1.3) node {$(\beta_2,\alpha_2)$};
\draw[color=black] (-3.8,1.3) -- (-3.8,0.7);
\fill [color=black] (-3.8,0.7) circle (1.5pt);
\draw[color=black] (-3.1,0.5) node {$(\beta_3,\alpha_3)$};
\draw[dashed,color=black] (-3.8,0.7) -- (-4.3,0.2);
\fill [color=black] (-4.3,0.2) circle (1.5pt);
\draw[color=black] (-5.2,0.3) node {$(\beta_n,\alpha_n)$};
\end{tikzpicture}

\caption{The sequence $(\beta_n,\alpha_n)$.}\label{fig:SequenceEn}

\end{center}
\end{figure}

We let $\widetilde{\alpha}\geq0$ be the limit of the sequence $\{\alpha_n\}$.
By assumption, we would get a contradiction if we prove that $\widetilde{\alpha}=0$.
Hence, assume this is not the case, namely $\widetilde{\alpha}>0$.
The idea is to find bounds for $\ch_0(E_n)$, $H^2\ch_1(E_n)$, and $H\ch_2(E_n)$.

\begin{Claim}\label{claim2}
For all $n>0$, the following inequality holds:
\[
\overline{\Delta}_H(E_n) + (\alpha_n H^3 \ch_0(E_n))^2 < \overline{\Delta}_H(E_0) + (\alpha_0 H^3 \ch_0(E_0))^2.
\]
\end{Claim}

\begin{Prf} (Claim \ref{claim2})
Again, we proceed by induction.
By Claim \ref{claim1}, and by definition of the generalized discriminant, we have
\begin{small}
\begin{equation*}
\begin{split}
\overline{\Delta}_H&(E_{n+1}) + (\alpha_{n+1} H^3 \ch_0(E_{n+1}))^2\\
    & = (H^2 \ch^{B_{n+1}}_1(E_{n+1}))^2 -2H^3 \ch_0 (E_{n+1}) (H\ch^{B_{n+1}}_2(E_{n+1})) + (\alpha_{n+1} H^3 \ch_0(E_{n+1}))^2\\
    & = (H^2 \ch^{B_{n+1}}_1(E_{n+1}))^2 -2H^3 \ch_0 (E_{n+1}) (\frac 12 \alpha_{n+1}^2 H^3 \ch_0 (E_{n+1})) +(\alpha_{n+1} H^3 \ch_0(E_{n+1}))^2\\
    & = (H^2 \ch^{B_{n+1}}_1(E_{n+1}))^2\\
    & < (H^2 \ch^{B_{n+1}}_1(E_n))^2\\
    & = (H^2 \ch^{B_{n+1}}_1(E_{n}))^2 -2H^3 \ch_0 (E_{n}) (\frac 12 \alpha_{n+1}^2 H^3 \ch_0 (E_{n})) +(\alpha_{n+1} H^3 \ch_0(E_{n}))^2\\
    & = (H^2 \ch^{B_{n+1}}_1(E_{n}))^2 -2H^3 \ch_0 (E_{n}) (H\ch^{B_{n+1}}_2(E_{n})) + (\alpha_{n+1} H^3 \ch_0(E_{n}))^2\\
    & = \overline{\Delta}_H(E_{n}) + (\alpha_{n+1} H^3 \ch_0(E_{n}))^2\\
    & \leq \overline{\Delta}_H(E_{n}) + (\alpha_{n} H^3 \ch_0(E_{n}))^2.
\end{split}
\end{equation*}
\end{small}
\end{Prf}

By Claim \ref{claim2}, we deduce, for all $n>0$, the inequality
\[
\overline{\Delta}_H(E_n) + (\widetilde{\alpha} H^3 \ch_0(E_n))^2 < \overline{\Delta}_H(E_0) + (\alpha_0 H^3 \ch_0(E_0))^2.
\]
Hence, we get immediately
\begin{equation}\label{eq:BoundDelta}
\overline{\Delta}_H(E_n) < \overline{\Delta}_H(E_0) + (\alpha_0 H^3 \ch_0(E_0))^2 =: \Gamma_0.
\end{equation}
and, by Theorem \ref{thm:BogNoCh2}, we have
\begin{equation}\label{eq:BoundRank}
(\ch_0(E_n))^2 < \frac{1}{(\widetilde{\alpha}H^3)^2} \left(\overline{\Delta}_H(E_0) + (\alpha_0 H^3 \ch_0(E_0))^2\right) = \Gamma_1.
\end{equation}
Finally, to bound $H^2\ch_1$, assume first that $\ch_0(E_n)\neq0$.
Then, by \eqref{eq:general}, \eqref{eq:BoundDelta}, \eqref{eq:BoundRank}, and Claim \ref{claim1}, we have
\begin{equation}\label{eq:BoundC1}
|H^2\ch_1(E_n)| \leq H^3 \sqrt{\Gamma_1}\left( |\beta_0| + |\alpha_0| + \sqrt{ \alpha_0^2 + \frac{\Gamma_0}{(H^3)^2} } \right) =: \Gamma_2.
\end{equation}
The case in which $\ch_0(E_n)=0$, follows by Claim \ref{claim1} by observing that, either $\ch_0(E_m)=0$, for all $0\leq m\leq n$, or there exists a maximum $0\leq m < n$ for which $\ch_0(E_m)\neq0$.
In the first case, we have
\begin{equation}\label{eq:BoundC1rk0Always}
0< H^2\ch_1(E_n) < H^2\ch_1(E_0),
\end{equation}
while in the second
\begin{equation}\label{eq:BoundC1rk0}
0 < H^2\ch_1(E_n) < |H^2\ch_1(E_m)| + |\beta_m| |\ch_0(E_m)| \leq \Gamma_2 + \left( |\beta_0| + |\alpha_0| \right) \Gamma_1.
\end{equation}

Summing up, by \eqref{eq:BoundDelta}, \eqref{eq:BoundRank}, \eqref{eq:BoundC1}, \eqref{eq:BoundC1rk0Always}, and \eqref{eq:BoundC1rk0}, we found bounds for $\ch_0(E_n)$, $H^2\ch_1(E_n)$, and $H\ch_2(E_n)$, for all $n$.
But this shows that these classes are finite, and so there must exist an object $E$ which does not satisfy the inequality in Conjecture \ref{con:strongBG} for all $\alpha$ close to $0$, which contradicts our assumption.
\end{Prf}

\section{The case of the projective space}\label{sec:P3}

In this section we expand \cite[Section 8.2]{BMT:3folds-BG} to show that in the case of $X=\P^3$, the assumptions in Proposition \ref{prop:reduction} are satisfied.
This will complete the proof of Theorem \ref{thm:main}.
To simplify notation, we directly identify $\Num_\R(\P^3)$ with $\R^{\oplus 4}$, and we take $\omega=\alpha,B=\beta\in\R$, $\alpha>0$.
The tilted slope becomes:
\[
\nu_{\alpha,\beta}=\frac{\ch_2^\beta-\frac{\alpha^2}{2} \ch_0}{\ch_1^\beta} = \frac{\ch_2-\beta\ch_1+\left(\frac{\beta^2}{2}-\frac{\alpha^2}{2}\right)\ch_0}{\ch_1-\beta \ch_0}.
\]

\begin{Prop}\label{prop:P3Bound}
For all $\alpha<\frac 13$, and for all $\beta\in\R$, Conjecture \ref{con:strongBG} holds.
\end{Prop}

The proof is an adaptation of \cite[Section 8.2]{BMT:3folds-BG}, where only the case $\beta=0$ was considered.
The idea is to use the existence of Bridgeland's stability conditions on $\Db(\P^3)$ associated to strong exceptional collections of sheaves (see \cite[Example 5.5]{Bridgeland:Stab} and \cite[Section 3.3]{Macri:curves}). Here, we will use the full strong exceptional collection $\EEE$ on $\Db(\P^3)$ given by
\[
\EEE := \left\{ \OO_{\P^3}(-1), \QQ, \OO_{\P^3}, \OO_{\P^3}(1)\right\},
\]
where $\QQ:=T_{\P^3}(-2)$ is given by
\[
0 \to \OO_{\P^3}(-2) \to \OO_{\P^3}(-1)^{\oplus 4} \to \QQ \to 0.
\]

We consider the region $V$ given by
\[
V:=\left\{ (\beta,\alpha)\in\mathbb{H}\,:\, \begin{array}{l}0\geq\beta>-\frac 23\\ 0<\alpha< \frac 13\end{array}\right\}.
\]

\begin{Lem}\label{lem:TensorProd}
Assume that Conjecture \ref{con:strongBG} holds, for all $(\beta,\alpha)\in V$.
Then it holds for all $\alpha<\frac 13$ and for all $\beta\in\R$.
\end{Lem}

\begin{Prf}
Assume, for a contradiction, there exist $\alpha_0<\frac 13$ and $\beta_0\in\R$, and $E\in\Db(\P^3)$ which does not satisfy Conjecture \ref{con:strongBG}.
By acting with the autoequivalence $\otimes\OO_{\P^3}(1)$ and with the local dualizing functor $\D(\blank) := \RlHom(\blank, \OO_X[1])$, we can assume (see \cite[Proposition 5.1.3]{BMT:3folds-BG}) that $0>\beta_0\geq-\frac 12$, which contradicts our assumption.
\end{Prf}

The next result will allow us to use the exceptional collection $\EEE$ for doing computations.
We postpone the proof to the end of the section.

\begin{Lem}\label{lem:Qstable}
For all $(\beta,\alpha)\in V$, we have $\QQ[1]\in\Coh^{\beta}(\P^3)$ and $\nu_{\alpha,\beta}^{\mathrm{min}}(\QQ[1])>0$.
\end{Lem}

We divide the region $V$ in three parts:
\begin{equation*}
\begin{split}
V_1 &:= \left\{ (\beta,\alpha)\in V\,:\,\beta < - \alpha \right\}\\
V_2 &:= \left\{ (\beta,\alpha)\in V\,:\,\beta > - \alpha \right\}\\
V_3 &:= \left\{ (\beta,\alpha)\in V\,:\,\beta = - \alpha \right\}.
\end{split}
\end{equation*}

\begin{figure}[htb]
\begin{center}
\includegraphics{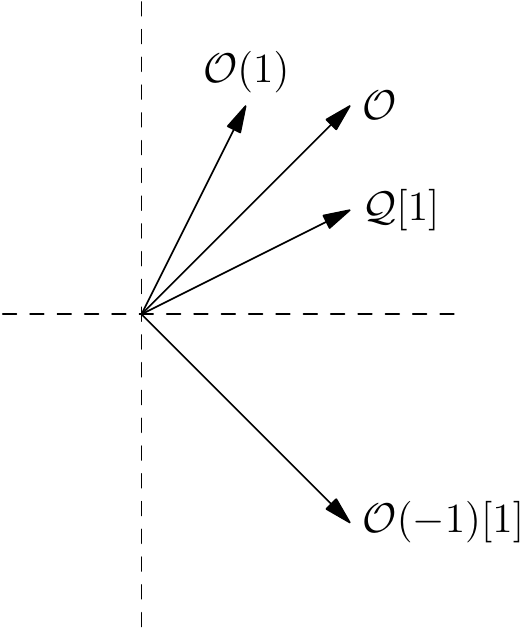}
\hspace{2cm}
\includegraphics{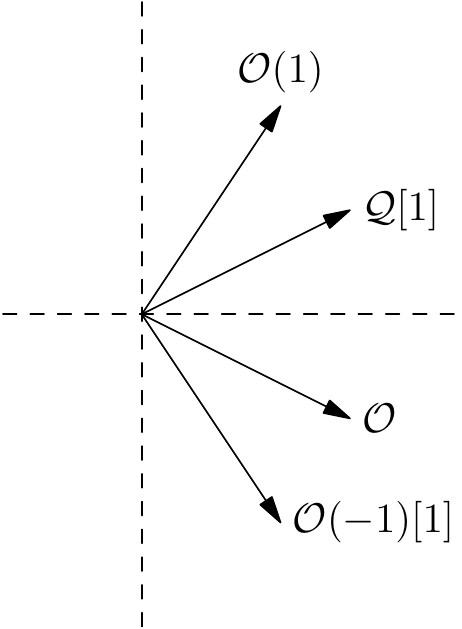}
\caption{The slopes in $\Coh^{\beta}(\P^3)$ of the exceptional objects when $(\beta,\alpha)\in V_1$ (left) and $(\beta,\alpha)\in V_2$ (right).
The tilt to $\AA^{\alpha,\beta}$ corresponds to considering the upper half-plane.}
\label{figure1}
\end{center}
\end{figure}

We first examine $V_1$ and $V_2$.
On $V_1$, we have:
\begin{equation*}
\begin{split}
&\nu_{\alpha,\beta}(\OO) = \frac 12\cdot \frac{\beta^2-\alpha^2}{-\beta}>0,\\
&\nu_{\alpha,\beta}(\OO(-1)) = \frac 12 \cdot \frac{(\beta+1)^2-\alpha^2}{-\beta-1}<0,\\
&\nu_{\alpha,\beta}(\OO(1)) = \frac 12\cdot \frac{(\beta-1)^2-\alpha^2}{1-\beta}>0,\\
&\nu_{\alpha,\beta}(\QQ) = \frac 32\cdot \frac{\left(\beta+\frac 23\right)^2-\alpha^2 - \frac 49}{-2-3\beta}>0.
\end{split}
\end{equation*}
On $V_2$ we get the same expressions, but now $\nu_{\alpha,\beta}(\OO)<0$ (see Figure \ref{figure1}).

We now tilt one more time $\Coh^{\beta}(\P^3)$, as explained in \cite[Definition 3.2.5]{BMT:3folds-BG}.
As in Section \ref{subsec:TiltStability}, we can define a torsion pair:
\begin{align*}
\TT_{\omega, B}' &= \stv{E \in \Coh^{\beta}(\P^3)}
{\text{any quotient $E \onto G$ satisfies $\nu_{\alpha,\beta}(G) > 0$}} \\
\FF_{\omega, B}' &= \stv{E \in \Coh^{\beta}(\P^3)}
{\text{any subsheaf $F \into E$ satisfies $\nu_{\alpha,\beta}(F) \le 0$}}.
\end{align*}
We let $\AA^{\alpha,\beta} \subset \Db(\P^3)$ be the extension-closure
\[
\AA^{\alpha,\beta} := \langle \TT_{\alpha,\beta}', \FF_{\alpha,\beta}'[1] \rangle.
\]
By the previous computation and Lemma \ref{lem:Qstable}, we have
\begin{align*}
& \left\{  \OO_{\P^3}(-1)[2], \QQ[1], \OO_{\P^3}, \OO_{\P^3}(1) \right\}\subset \AA^{\alpha,\beta}, & \text{ for }(\beta,\alpha)\in V_1,\\
& \left\{  \OO_{\P^3}(-1)[2], \QQ[1], \OO_{\P^3}[1], \OO_{\P^3}(1) \right\}\subset \AA^{\alpha,\beta}, & \text{ for }(\beta,\alpha)\in V_2.\\
\end{align*}

On the category $\AA^{\alpha,\beta}$ we consider the following function (a posteriori, this will be a slope function):
\[
\lambda_{\alpha,\beta}:=\begin{cases} 
 + \infty, & \text{ if }  \ch_2^\beta-\frac{\alpha^2}{2} \ch_0^\beta= 0,\\
\ & \\
\frac{\ch_3^\beta-\frac{\alpha^2}{6}\ch_1^\beta}{\ch_2^\beta-\frac{\alpha^2}{2} \ch_0^\beta}, & \text{ otherwise}.
\end{cases}
\]

We have:
\begin{equation*}
\begin{split}
&\lambda_{\alpha,\beta}(\OO) = -\frac{\beta}{3},\\
&\lambda_{\alpha,\beta}(\OO(-1)) = -\frac{\beta}{3}-\frac 13,\\
&\lambda_{\alpha,\beta}(\OO(1)) = -\frac{\beta}{3} + \frac 13,\\
&\lambda_{\alpha,\beta}(\QQ) = \frac{\left( \frac 23 - \beta^2 - \frac{\beta^3}{2}\right) + \frac{\alpha^2}{6} \left( 3\beta + 2\right)}{2\beta + \frac 32 \beta^2 - \frac 32 \alpha^2}.
\end{split}
\end{equation*}

\begin{figure}[htb]
\begin{center}
\includegraphics{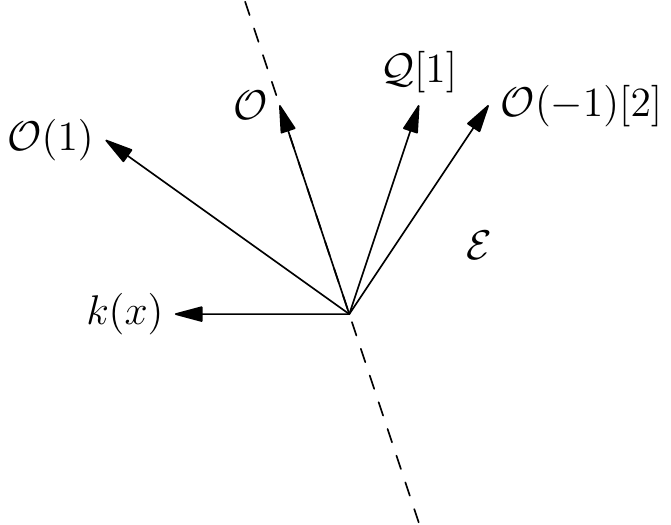}
\hspace{1cm}
\includegraphics{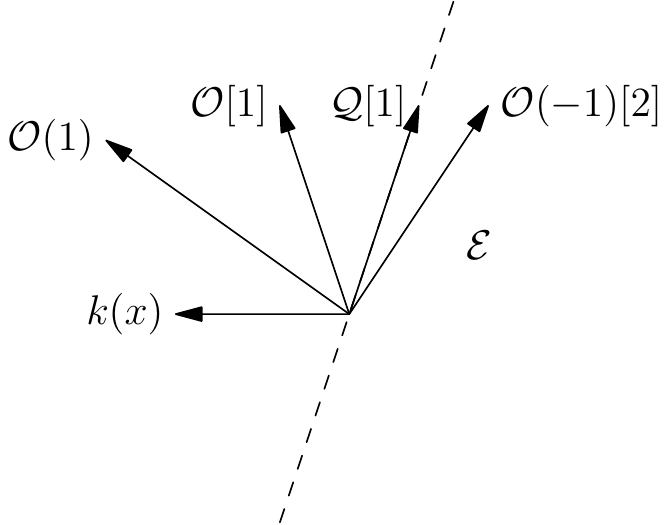}
\caption{The slopes in $\AA^{\alpha,\beta}$ of the exceptional objects and the skyscraper sheaves when $(\beta,\alpha)\in V_1$ (left) and $(\beta,\alpha)\in V_2$ (right).
The category $\EE$, obtained by tilting to the right along the dotted line, is the extension-closed subcategory generated by $\OO(-1)[2]$, $\QQ[1]$, $\OO$, $\OO(1)[-1]$.
It is equivalent to the category of modules over the finite-dimensional algebra determined by the dual exceptional collection to $\EEE$.}
\label{figure2}
\end{center}
\end{figure}

On $V_1$, we deduce that $\lambda_{\alpha,\beta} (Q) < \lambda_{\alpha,\beta} (\OO(1))$, while, on $V_2$, $\lambda_{\alpha,\beta} (Q) < \lambda_{\alpha,\beta} (\OO)$ (see Figure \ref{figure2}).

By \cite[Proposition 8.1.1]{BMT:3folds-BG} (and mimicking the proof of \cite[Theorem 8.2.1]{BMT:3folds-BG}), this shows that Conjecture \ref{con:strongBG} holds, for all $(\beta,\alpha)\in V_1 \cup V_2$.

To deal with the region $V_3$ (namely, the case $\alpha=-\beta$), we consider a slightly modified function on $\AA^{\alpha,\beta}$:
\[
\lambda_{\alpha,\beta}:=\begin{cases} 
 + \infty, & \text{ if }  \ch_2^\beta-\frac{\alpha^2}{2} \ch_0^\beta= 0,\\
\ & \\
\frac{\ch_3^\beta-\frac{\alpha^2}{6}\ch_1^\beta - \epsilon\ch_1^\beta}{\ch_2^\beta-\frac{\alpha^2}{2} \ch_0^\beta}, & \text{ otherwise},
\end{cases}
\]
where $\epsilon>0$.
In this case, we still have
\[
\left\{  \OO_{\P^3}(-1)[2], \QQ[1], \OO_{\P^3}[1], \OO_{\P^3}(1) \right\}\subset \AA^{\alpha,\beta},
\]
and
\begin{equation*}
\begin{split}
&\lambda_{\alpha,\beta}(\OO) = +\infty,\\
&\lambda_{\alpha,\beta}(\OO(-1)) = -\frac{\beta}{3}-\frac 13 + 2\epsilon \cdot \frac{\beta+1}{2\beta+1},\\
&\lambda_{\alpha,\beta}(\OO(1)) = -\frac{\beta}{3} + \frac 13 + 2\epsilon \cdot \frac{\beta-1}{1-2\beta},\\
&\lambda_{\alpha,\beta}(\QQ) = \frac{1-\beta^2}{3\beta} + \epsilon \cdot \frac{3\beta+2}{2\beta}.
\end{split}
\end{equation*}
We deduce that, for all $0>\beta>-\frac 13$, there exists $\epsilon(\beta)>0$ such that
\begin{align*}
\lambda_{\alpha,\beta}(\OO(1)) > \lambda_{\alpha,\beta}(\OO(-1)) \quad \text{ and } \quad \lambda_{\alpha,\beta}(\OO(1)) > \lambda_{\alpha,\beta}(\QQ),
\end{align*}
when $(\beta,\alpha)\in V_3$ and $\epsilon<\epsilon(\beta)$.
Again, by \cite[Proposition 8.1.1]{BMT:3folds-BG}, if we fix $\beta$ and let $\epsilon\to 0$, this shows that Conjecture \ref{con:strongBG} holds also for all $(\beta,\alpha)\in V_3$.
By Lemma \ref{lem:TensorProd}, this would complete the proof of Proposition \ref{prop:P3Bound}, once Lemma \ref{lem:Qstable} is proved.

\begin{Prf} (Lemma \ref{lem:Qstable})
Since $\QQ\in\Coh(\P^3)$ is slope-stable, with Chern character $\ch(\QQ)=(3,-2,0,\frac 23)$, we have, by definition, $\QQ[1]\in\Coh^{\beta}(\P^3)$, for all $\beta\geq-\frac 23$.
Moreover, for $0\geq\beta >-\frac 23$ and for all $\alpha>0$, we have $\nu_{\alpha,\beta}(\QQ[1])>0$.

Assume, for a contradiction, that there exists $(\beta_0,\alpha_0)\in V$ such that $\nu_{\alpha_0,\beta_0}^{\mathrm{min}}(Q[1])\leq0$.
Let $N_0\in\Coh^{\beta_0}(\P^3)$ be the tilt-stable quotient $\QQ[1]\onto N_0$ in $\Coh^{\beta_0}(\P^3)$ such that $\nu_{\alpha_0,\beta_0}(N_0)\leq0$.
By taking the long exact sequence in cohomology, $N_0\cong M_0[1]$, where $M_0\in\Coh(\P^3)$ is a torsion-free sheaf.

Consider the curves $\CC_0$, given by $\nu_{\alpha,\beta}(N_0)=0$, in the region $\beta>\frac{\ch_1(M_0)}{\ch_0(M_0)}$, and $\BB_0$, given by $\nu_{\alpha,\beta}(\QQ[1])=\nu_{\alpha,\beta}(N_0)$.
Since the vector $(3,-2,0)$ is primitive, $\BB_0$ must be a semi-circle in $\H$.
Consider the unique point of intersection $(x,y)\in \CC_0\cap \BB_0$.
Since $\nu_{\alpha,\beta}(\QQ[1])>0$, for $0\geq\beta >-\frac 23$, we have $x>0$.
In particular, $\BB_0\cap\{\beta=0 \}\neq\emptyset$.
See Figure \ref{fig:StabilityQ}.

\begin{figure}
\begin{center}

\definecolor{uququq}{rgb}{0.25,0.25,0.25}
\begin{tikzpicture}[line cap=round,line join=round,>=triangle 45,x=1.0cm,y=1.0cm]
\draw[->,color=black] (-5,0) -- (5,0);
\draw[->,color=black] (0,-0.5) -- (0,5);
\fill [color=black] (0,0) circle (1.5pt);
\draw[color=black] (0.2,-0.2) node {$0$};
\draw[color=black] (5.1,-0.2) node {$\beta$};
\draw[color=black] (0.2,5.1) node {$\alpha$};
\clip(-3.5,0) rectangle (5.5,4);
\draw (2.2,3.8) node[anchor=north west] {$\CC_0$};
\draw (-2.2,3) node[anchor=north west] {$\BB_0$};
\draw [samples=50,domain=-0.99:0.99,xshift=0cm,yshift=0cm] plot
({-3+1.73*(1+(\x)^2)/(1-(\x)^2)},{2*\x/(1-\x^2)});
\draw [samples=50,domain=-2:2,xshift=0cm,yshift=0cm] plot
({-0.5+2.5*((2*((0.5)*\x))/(1+((0.5)*\x)^2))},{2.5*((1-((0.5)*\x)^2)/(1+((0.5)*\x)^2))});
\fill [color=black] (0.95,2.05) circle (1.5pt);
\draw[color=black] (1.8,2.05) node {$(x,y)$};
\end{tikzpicture}

\caption{The curves $\BB_0$ and $\CC_0$.}\label{fig:StabilityQ}

\end{center}
\end{figure}
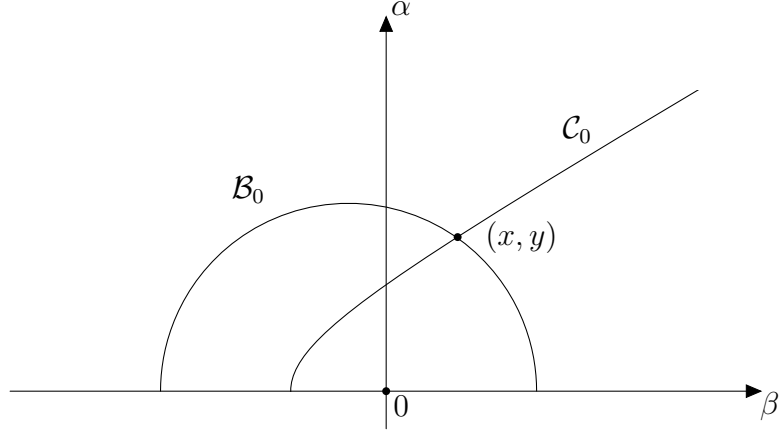

By Bertram's Nested Wall Theorem of \cite{Maciocia:walls} (whose proof works as well in our context, due to Theorem \ref{thm:BogNoCh2}), we know that \emph{pseusdo-walls} for $\QQ[1]$ are \emph{nested} semi-circles, namely, either $\QQ[1]$ is tilt-stable outside $\BB_0$ and unstable in the interior, or there exists another semi-circle $\BB_1$ with the same property and $\BB_1$ contains $\BB_0$ in its interior.
In both cases, by the previous argument, the semi-circles $\BB_0$ and $\BB_1$ intersect the semi-line $\beta=0$.
Hence, there exists $\alpha_1>0$ such that $\QQ[1]$ is not $\nu_{\alpha_1,0}$-stable.
This contradicts Lemma \ref{lem:auxil} below.
\end{Prf}

\begin{Lem}\label{lem:auxil}
For all $\alpha>0$, $\QQ[1]$ is $\nu_{\alpha,0}$-stable.
\end{Lem}

\begin{Prf}
First of all, we observe that $\QQ[1]$ is $\mathrm{PGL}(4)$-invariant.
By uniqueness of Harder-Narasimhan filtrations, if $\QQ[1]$ is not tilt-stable, then its HN factors have to be $\mathrm{PGL}(4)$-invariant as well.

Consider the category $\Coh^{\beta=0}(\P^3)$.
The function $f_{0}:=\ch_1$ is additive and takes positive integral values on $\Coh^{0}(\P^3)$.
Since $f_0(\QQ[1])=2$, if there exists an exact sequence in $\Coh^{0}(\P^3)$
\begin{equation}\label{eq:destab}
0\to P \to \QQ[1] \to N\cong M[1] \to 0
\end{equation}
which is destabilizing, then $f_0(P)=f_0(N)=1$ and both $P$ and $N$ must be tilt-stable.
The long exact sequence in cohomology gives
\[
0\to \HH^{-1}(P)\to\QQ\to M\to\HH^0(P)\to 0,
\]
with $\HH^{-1}(P)$ and $M$ torsion-free with $\mu^{\mathrm{max}}_{\alpha,0}\leq0$.
Since \eqref{eq:destab} is destabilizing, we must have $\mu^{\mathrm{max}}_{\alpha,0}(M),\mu^{\mathrm{max}}_{\alpha,0}(\HH^{-1}(P))<0$.
This shows that there are only two possibilities:
\begin{enumerate}
\item \label{case:1} either $\ch_1(M)=\ch_1(\HH^{-1}(P))=-1$,
\item \label{case:2} or $\HH^{-1}(P)=0$.
\end{enumerate}

For Case \eqref{case:1}, we must have $\ch_1(\HH^0(P))=0$, and so $\HH^0(P)$ is a torsion sheaf supported on a one-dimensional subscheme.
By the $\mathrm{PGL}(4)$-invariance, $\HH^0(P)=0$.
Finally, since $\QQ$ is slope-stable, we must have $\ch_0(\HH^{-1}(P))=1$, and so $\HH^{-1}(P)\cong\II_C(-1)$, for $C\subset\P^3$ a one-dimensional subscheme of degree $d\geq0$.
Since $\Hom(\OO_C(-1)[-1],\QQ)=0$, then the inclusion $\II_C(-1)\into\QQ$ factorizes through $\OO(-1)\into\QQ$ and so, by the Snake Lemma, $M$ has torsion, which is a contradiction, unless $C=0$.
Summarizing, we proved that in Case \eqref{case:1}, $P\cong\OO_{\P^3}(-1)[1]$.
But then, the equation $\nu_{\alpha,0}(\QQ[1])=\nu_{\alpha,0}(P)$ has no solutions, and so \eqref{eq:destab} cannot be destabilizing.

For Case \eqref{case:2}, we have $P\in\Coh(\P^3)$ and an exact sequence in $\Coh(\P^3)$
\[
0\to \QQ \to M \to P\to 0,
\]
with $\ch_1(M)=-1$, $\ch_1(P)=1$, and $\ch_0(M)\geq3$.
We now use Theorem \ref{thm:BogNoCh2} once more.
Indeed, since $N$ must be tilt-stable, we have
\[
\ch_2(M)\leq \frac{1}{2\ch_0(M)},
\]
and so $\ch_2(M)\leq0$.
As a consequence, the equation $\nu_{\alpha,0}(\QQ[1])=\nu_{\alpha,0}(P)$ has no solutions $\alpha>0$, and so \eqref{eq:destab} cannot be destabilizing also in this case.
\end{Prf}

\section{An application}\label{sec:applications}

In this section we briefly discuss an application of Theorem \ref{thm:main}, and some examples.

In \cite[Example 7.2.4]{BMT:3folds-BG}, we pointed out a relation between Conjecture \ref{con:strongBG} and Castelnuovo's inequality for curves in $\P^3$.
In particular, by using Castelnuovo's inequality, we showed that Conjecture \ref{con:strongBG} holds for ideal sheaves of curves with respect to some tilt-stability.
It is interesting to observe that a sort of vice versa holds: from Theorem \ref{thm:main}, we can deduce a certain inequality for curves in $\P^3$, which is much weaker than Castelnuovo's one, but already non-trivial:

\begin{Cor}\label{cor:Castelnuovo}
Let $C$ be a pure one-dimensional scheme in $\P^3$ of degree $d\geq2$.
Let $h:=\ch_3(\II_C) - 2d$.
Then
\begin{equation}\label{eq:Castelnuovo1}
h \leq \frac{2d^2-5d}{3}.
\end{equation}
Moreover, if $C$ is integral and not contained in a plane, then
\begin{equation}\label{eq:Castelnuovo2}
h \leq \frac{d^2-4d}{3}.
\end{equation}
\end{Cor}

We recall that, for an ideal sheaf $\II_C$ of an integral curve $C\subset\P^3$ of degree d and arithmetic genus $g$, $h=g-1$.
Hence, the inequality \eqref{eq:Castelnuovo2} compares with \cite[IV, 6.4]{Hartshorne}.

To prove Corollary \ref{cor:Castelnuovo}, we introduce some more notation.
We denote by $\BB_1$ and $\BB_2$ the two semi-circles
\begin{align*}
\BB_1 &:\quad \alpha^2 +\left( \beta + \frac{2d+1}{2}\right)^2 = \left(\frac{2d-1}{2}\right)^2\\
\BB_2 &:\quad \alpha^2 +\left( \beta + \frac{d+2}{2}\right)^2 = \left(\frac{d-2}{2}\right)^2.
\end{align*}
They correspond to the loci $\nu_{\alpha,\beta}(\II_C)=\nu_{\alpha,\beta}(\OO_{\P^3}(-1))$ and $\nu_{\alpha,\beta}(\II_C)=\nu_{\alpha,\beta}(\OO_{\P^3}(-2))$, respectively.
More generally, for an object $A\in\Db(\P^3)$ such that $(\ch_0(A),\ch_1(A),\ch_2(A))$ is not a multiple of $(1,0,-d)$, we denote by $\BB_A$ the semi-circle with equation $\nu_{\alpha,\beta}(\II_C)=\nu_{\alpha,\beta}(A)$.

Finally, as in Section \ref{subsec:bound}, we denote by $\CC$ the branch of the hyperbola $\nu_{\alpha,\beta}(\II_C)=0$ in $\H$; explicitly,
\[
\CC:\quad \beta^2 - \alpha^2 = 2d,\qquad \beta<0.
\]

\begin{Prf} (Corollary \ref{cor:Castelnuovo})
For the first part of the statement, we would like to show that on the exterior part of the semi-circle $\BB_1$ in $\H\cap\{-2d<\beta<-1\}$ the ideal sheaf $\II_C$ is $\nu_{\alpha,\beta}$-stable.

First of all, we consider the semi-line $\beta=-1$ and the category $\Coh^{\beta=-1}(\P^3)$.
The function $f_{-1}:=\ch_1 + \ch_0$ is additive and takes positive integral values on $\Coh^{-1}(\P^3)$.
Since $f_{-1}(\II_C)=1$, then $\II_C$ must be $\nu_{\alpha,-1}$-stable, for all $\alpha>0$.

We now consider the semi-line $\beta=-2$ and the category $\Coh^{-2}(\P^3)$.
By \cite[Proposition 14.2]{Bridgeland:K3} (whose proof generalizes to our case), we know that, for $\alpha\gg0$, $\II_C$ is $\nu_{\alpha,-2}$-stable.
Assume that $\II_C$ is not $\nu_{\alpha,-2}$-semistable, for all $\alpha>0$.
Then, by Proposition \ref{prop:openness}, there exists $\alpha_0>0$ such that $\II_C$ is $\nu_{\alpha,-2}$-stable, for $\alpha>\alpha_0$, it is $\nu_{\alpha,-2}$-semistable at $\alpha=\alpha_0$, and not semistable for $\alpha<\alpha_0$.
Then $\alpha_0$ must be in the intersection of the semi-line $\beta=-2$ with a semi-circle $\BB_A$, for some $A\in\Coh^{-2}(\P^3)$ such that $A\into\II_C$ in $\Coh^{-2}(\P^3)$.
By looking at the long exact sequence in cohomology, we deduce that $A\in\Coh(\P^3)$, $\ch_0(A)\geq1$, and it is torsion-free.
Moreover, since the function $f_{-2}:=\ch_1 + 2\ch_0$ is additive, takes positive integral values on $\Coh^{-2}(\P^3)$, and $f_{-2}(\II_C)=2$, we must have $f_{-2}(A)=1$, namely
\[
\frac{\ch_1(A)}{\ch_0(A)} = -2 + \frac{1}{\ch_0(A)}.
\]

Let $(-2,\alpha_1)$ be the intersection point in $\H$ between $\beta=-2$ and $\BB_1$ (the intersection is non-empty, since $d\geq2$).
We claim that $\alpha_0\leq\alpha_1$.
Indeed, if $\ch_0(A)=1$, then $\ch_1(A)=-1$.
Hence, $A\cong\II_W(-1)$, for some subscheme $W$ of dimension $1$.
Therefore, $\alpha_0\leq\alpha_1$.
If $\ch_0(A)\geq2$, then $-2<\frac{\ch_1(A)}{\ch_0(A)}<-1$.
By Bertram's Nested Wall Theorem of \cite{Maciocia:walls}, we know that either $\BB_A=\BB_1$, or they are disjoint.
Since $\BB_A\cap\{\beta=\frac{\ch_1(A)}{\ch_0(A)}\}=\emptyset$, this immediately implies that $\alpha_0\leq\alpha_1$, as we wanted.

By using the Nested Wall Theorem again, since we proved that, on the line $\beta=-1$, the ideal sheaf $\II_C$ is stable and, on the line $\beta=-2$, the first wall is $\BB_1$, this shows that on the exterior part of the semi-circle $\BB_1$ in $\H\cap\{-2d<\beta<-1\}$ the ideal sheaf $\II_C$ is $\nu_{\alpha,\beta}$-stable, which is what we wanted.
To get the inequality \eqref{eq:Castelnuovo1}, we only need to compute the intersection point $\CC\cap\BB_1$.
Theorem \ref{thm:main} yields then directly \eqref{eq:Castelnuovo1}.

The proof of \eqref{eq:Castelnuovo2} is very similar.
We consider the semi-line $\beta=-3$, the category $\AA_{-3}:=\Coh^{-3}(\P^3)$, and $A\into\II_C$ in $\Coh^{-3}(\P^3)$.
By looking at the function $f_{-3}:=\ch_1+3\ch_0$, we must have either $f_{-3}(A)=1$, or $=2$.
If $\ch_0(A)\geq3$, then by using again \cite{Maciocia:walls}, we can deduce that $\BB_A$ is contained in the interior of $\BB_2$.
If $\ch_0(A)=2$, we distinguish two possibilities, according whether $f_{-3}(A)=1$, or $=2$.
If $=1$, then we can argue as before, and deduce that $\BB_A$ is contained in the interior of $\BB_2$.
If $=2$, then $\ch_1(A)=-4$, and so, by Theorem \ref{thm:BogNoCh2}, $\ch_2(A)\leq4$.
If $\ch_2(A)=4$, then $\BB_A=\BB_2$.
If $\ch_2(A)<4$, then $\BB_A$ is again contained in the interior of $\BB_2$.

Finally, if $\ch_0(A)=1$, then either $A\cong\II_W(-2)$, or $A\cong\II_W(-1)$, with $W$ a closed subscheme of dimension $1$.
The first case, can be dealt as before.
To exclude the second case, we use the assumption that $C$ is integral and not contained in a plane.
Indeed, in such a case, we must have $C\subset W$, and so $A\into\II_C$ does not destabilize.

As before, to get the inequality \eqref{eq:Castelnuovo2}, we only need to compute the intersection point $\CC\cap\BB_2$ and apply Theorem \ref{thm:main}.
\end{Prf}

\begin{Ex}\label{ex:d=1}
For the case $d=1$, the situation is slightly degenerate.
Indeed, in such a case, $\II_C$ is $\nu_{\alpha,\beta}$-semistable for all $(\beta,\alpha)\in\H$ for which
\[
\alpha^2 + \left(\beta+\frac 32\right)^2 \geq \frac 14.
\]
Hence, in particular, it is semistable for all $(\beta,\alpha)\in\CC$.
Theorem \ref{thm:main} gives then $h\leq -\frac23$, namely $g (=0) \leq \frac 13$.
\end{Ex}

\begin{Ex}\label{ex:BeyondCastelnuovo}
The corresponding statement as Corollary \ref{cor:Castelnuovo} when the curve $C$ is contained in a surface $F\subset \P^3$ of degree $k>0$, and it is not contained in any surface of smaller degree, is not clear anymore.
In particular, it does not follow from the strong Castelnuovo's Theorem, proved by Harris \cite{Harris:space-curves,Hartshorne:CastelnuovoI}.

To be precise, it is not true that the first wall when $\II_C$ is destabilized coincide with the locus
\[
\nu_{\alpha,\beta}(\II_C) = \nu_{\alpha,\beta}(\OO_{\P^3}(-k)), \quad \text{ namely, } \quad \OO_{\P^3}(-k)\into\II_C.
\]
The simplest example (cfr.~\cite[V, 4.13.1]{Hartshorne}) is when $C$ is smooth with $k=3$, $d=7$, $g=5$.
In such a case, a destabilizing quotient is given instead by
\[
\II_C \onto \OO_{\P^3}(-5)[1].
\]
This gives the (well-known) existence of a non-trivial extension, $\GG\in\Coh(\P^3)$ of rank $2$, which must be stable.
It may be interesting to study the general situation, and see which kind of new stable objects arise as destabilizing factors of $\II_C$.
\end{Ex}

\bibliography{all}                      

\begin{thebibliography}{BBMT11}

\bibitem[ABCH12]{ABCH:MMP}
Daniele Arcara, Aaron Bertram, Izzet Coskun, and Jack Huizenga.
\newblock The minimal model program for {H}ilbert schemes of points on the
  projective plane and {B}ridgeland stability, 2012.
\newblock arXiv:1203.0316.

\bibitem[ABL07]{Aaron-Daniele}
Daniele Arcara, Aaron Bertram, and Max Lieblich.
\newblock Bridgeland-stable moduli spaces for {K}-trivial surfaces, 2007.
\newblock arXiv:0708.2247.

\bibitem[Bay10]{stability-tour}
Arend Bayer.
\newblock A tour to stability conditions on derived categories, 2010.
\newblock Available at {\tt http://www.math.uconn.edu/$\sim$bayer/}.

\bibitem[BBMT11]{BBMT:Fujita}
Arend Bayer, Aaron Bertram, Emanuele Macr{\`{\i}}, and Yukinobu Toda.
\newblock Bridgeland stability conditions on threefolds {II}: {A}n application
  to {F}ujita's conjecture, 2011.
\newblock arXiv:1106.3430.

\bibitem[BM12]{BM:Projectivity}
Arend Bayer and Emanuele Macr{\`{\i}}.
\newblock Projectivity and birational geometry of {B}ridgeland moduli spaces,
  2012.
\newblock arXiv:1203.4613.

\bibitem[BMT11]{BMT:3folds-BG}
Arend Bayer, Emanuele Macr{\`{\i}}, and Yukinobu Toda.
\newblock Bridgeland stability conditions on threefolds {I}:
  {B}ogomolov-{G}ieseker type inequalities, 2011.
\newblock To appear in \emph{J. Alg. Geom.}
\newblock arXiv:1103.5010.

\bibitem[Bri07]{Bridgeland:Stab}
Tom Bridgeland.
\newblock Stability conditions on triangulated categories.
\newblock {\em Ann. of Math. (2)}, 166(2):317--345, 2007.
\newblock arXiv:math/0212237.

\bibitem[Bri08]{Bridgeland:K3}
Tom Bridgeland.
\newblock Stability conditions on {$K3$} surfaces.
\newblock {\em Duke Math. J.}, 141(2):241--291, 2008.
\newblock arXiv:math/0307164.

\bibitem[Bri09]{Bridgeland:spaces}
Tom Bridgeland.
\newblock Spaces of stability conditions.
\newblock In {\em Algebraic geometry---{S}eattle 2005. {P}art 1}, volume~80 of
  {\em Proc. Sympos. Pure Math.}, pages 1--21. Amer. Math. Soc., Providence,
  RI, 2009.
\newblock arXiv:math.AG/0611510.

\bibitem[Har77]{Hartshorne}
Robin Hartshorne.
\newblock {\em Algebraic geometry}.
\newblock Springer-Verlag, New York, 1977.
\newblock Graduate Texts in Mathematics, No. 52.

\bibitem[Har78]{Hartshorne:CastelnuovoI}
Robin Hartshorne.
\newblock Stable vector bundles of rank {$2$} on {${\bf P}^{3}$}.
\newblock {\em Math. Ann.}, 238(3):229--280, 1978.

\bibitem[Har80]{Harris:space-curves}
Joe Harris.
\newblock The genus of space curves.
\newblock {\em Math. Ann.}, 249(3):191--204, 1980.

\bibitem[HRS96]{Happel-al:tilting}
Dieter Happel, Idun Reiten, and Sverre~O. Smal{\o}.
\newblock Tilting in abelian categories and quasitilted algebras.
\newblock {\em Mem. Amer. Math. Soc.}, 120(575):viii+ 88, 1996.

\bibitem[Huy12]{Huybrechts:Survey}
Daniel Huybrechts.
\newblock Introduction to stability conditions, 2012.
\newblock arXiv:1111.1745.

\bibitem[LQ11]{LoQin:miniwalls}
Jason Lo and Zhenbo Qin.
\newblock Mini-walls for {B}ridgeland stability conditions on the derived
  category of sheaves over surfaces, 2011.
\newblock arXiv:1103.4352.

\bibitem[Mac07]{Macri:curves}
Emanuele Macr{\`{\i}}.
\newblock Stability conditions on curves.
\newblock {\em Math. Res. Lett.}, 14(4):657--672, 2007.
\newblock arXiv:0705.3794.

\bibitem[Mac12]{Maciocia:walls}
Antony Maciocia.
\newblock Computing the walls associated to {B}ridgeland stability conditions
  on projective surfaces, 2012.
\newblock arXiv:1202.4587.

\bibitem[MM11]{MaciociaMeachan}
Antony Maciocia and Ciaran Meachan.
\newblock Rank one {B}ridgeland stable moduli spaces on a principally polarized
  abelian surface, 2011.
\newblock arXiv:1107.5304.

\bibitem[MYY11]{MYY2}
Hiroki Minamide, Shintarou Yanagida, and K{\=o}ta Yoshioka.
\newblock Some moduli spaces of {B}ridgeland's stability conditions, 2011.
\newblock arXiv:1111.6187.

\bibitem[Pol12]{Polishchuk:LIobjects}
Alexander Polishchuk.
\newblock Phases of {L}agrangian-{I}nvariant objects in the derived category of
  an abelian variety, 2012.
\newblock arXiv:1203.2300.

\bibitem[Tod11a]{Toda:DenefMoore}
Yukinobu Toda.
\newblock {B}ogomolov-{G}ieseker type inequality and counting invariants, 2011.
\newblock arXiv:1112.3411.

\bibitem[Tod11b]{Toda:Survey}
Yukinobu Toda.
\newblock Introduction and open problems of {D}onaldson-{T}homas theory, 2011.
\newblock To appear in Proceedings of the Conference ``Derived Categories'',
  Tokyo 2011.

\bibitem[Tod12a]{Toda:MMPxSurfaces}
Yukinobu Toda.
\newblock Stability conditions and birational geometry of projective surfaces,
  2012.
\newblock arXiv:1205.3602.

\bibitem[Tod12b]{Toda:ExtremalContractions}
Yukinobu Toda.
\newblock Stability conditions and extremal contractions, 2012.
\newblock arXiv:1204.0602.

\bibitem[Yos12]{Yoshioka:GeneralizedKummer}
K{\=o}ta Yoshioka.
\newblock Bridgeland's stability and the positive cone of the moduli spaces of
  stable objects on an abelian surface, 2012.
\newblock arXiv:1206.4838.

\bibitem[YY12]{YY:abeliansurfaces}
Shintarou Yanagida and K{\=o}ta Yoshioka.
\newblock Bridgeland's stabilities on abelian surfaces, 2012.
\newblock arXiv:1203.0884.

\end{thebibliography}
\bibliographystyle{alphaspecial}     

\end{document}